\documentclass[10pt,a4paper]{article}

\usepackage{amsmath}
\usepackage{graphicx,latexsym,euscript,makeidx,color}
\usepackage{amsmath,amsfonts,amssymb,amsthm,mathrsfs}
\usepackage{enumerate}

\usepackage{graphicx,colordvi}

\usepackage{amssymb}
\usepackage{amsmath}
\usepackage{amsfonts}
\usepackage{graphicx}
\usepackage{graphicx}          
\usepackage{color}

\usepackage[colorlinks,linkcolor=blue,anchorcolor=green,citecolor=red]{hyperref}

\usepackage{geometry}
\font\tenbb=msbm10 \font\sevenbb=msbm7 \font\fivebb=msbm5
\newfam\bbfam
\scriptscriptfont\bbfam=\fivebb \textfont\bbfam=\tenbb
\scriptfont\bbfam=\sevenbb

\newtheorem{theorem}{\indent Theorem}[section]
\newtheorem{definition}[theorem]{\indent Definition}
\newtheorem{proposition}[theorem]{\indent Proposition}
\newtheorem{corollary}[theorem]{\indent Corollary}

\newtheorem{remark}[theorem]{\indent Remark}

\makeatletter
   
   \@addtoreset{equation}{section}
\makeatother

\allowdisplaybreaks[4]

\begin{document}

\title{\bf Equilibrium Solutions of Multi-Period Mean-Variance Portfolio Selection
\thanks{This work is supported in part by the National Natural Science Foundation of China (61227902, 11471242, 61773222), the National Key Basic Research Program (973 Program) of China (2014CB845301), Hong Kong RGC grants 15224215 and 15255416.}
}
\author{Yuan-Hua Ni\thanks{College of Computer and Control Engineering, Nankai University, Tianjin 300350, P.R. China. Email: {\tt yhni@nankai.edu.cn}.}~~~~~~Xun Li\thanks{Department of Applied Mathematics, The Hong Kong
Polytechnic University, Kowloon, Hong Kong, P.R. China. Email: {\tt malixun@polyu.edu.hk}.}~~~~~~Ji-Feng Zhang\thanks{Key Laboratory of Systems and Control,
Institute of Systems Science, Academy of Mathematics and Systems Science,
Chinese Academy of Sciences, Beijing 100190; School of
Mathematical Sciences, University of Chinese Academy of Sciences, Beijing
100049, P.R. China. Email: {\tt jif@iss.ac.cn}.}~~~~~~Miroslav Krstic\thanks{Department of Mechanical and Aerospace Engineering, University of California, San Diego, USA. Email: {\tt krstic@ucsd.edu}.}}
\maketitle

{\bf Abstract:} This is a companion paper of [Mixed equilibrium solution of time-inconsistent stochastic LQ problem, arXiv:1802.03032], where general theory has been  established to characterize the open-loop equilibrium control, feedback equilibrium strategy and mixed equilibrium solution for a time-inconsistent stochastic linear-quadratic problem.
This note is, on the one hand to test the developed theory of that paper, and on the other hand to push the solvability of multi-period mean-variance portfolio selection. A nondegenerate assumption has been removed in this note, which is popular in existing literature about multi-period mean-variance portfolio selection; and neat conditions have been obtained to characterize the existence of  equilibrium solutions.


\section{Introduction}

Recently, a notion named mixed equilibrium solution is introduced in \cite{Ni-Li-Zhang-Krstic2018} for the time-inconsistent stochastic linear-quadratic (LQ, for short) optimal control; it contains two different parts: a pure-feedback-strategy part and an open-loop-control part, which together constitute a time-consistent solution. It is shown that the open-loop-control part will be of the feedback form of the equilibrium state. If we let the pure-feedback-strategy part be zero or let the open-loop-control part be independent of the initial state, then the mixed equilibrium solution will reduce to the open-loop equilibrium control and the (linear) feedback equilibrium strategy, respectively, both of which have been extensively studied in existing literature \cite{Basak,{Bjork-2},Cui-2017,Cui-2012,Hu-jin-Zhou,Hu-jin-zhou-2,Ni-Zhang-Krstic2017,Qi-Zhang2017,Yong-1,Yong-0,Yong-2013}.
Furthermore, the mixed equilibrium solution is not a hollow concept, whose study will give us more flexibility  to deal with the time-inconsistent optimal control.

The multi-period mean-variance portfolio selection is a particular example of time-inconsistent problem. In fact, the recent developments in time-inconsistent problems and the revisits of multi-period mean-variance portfolio selection \cite{Basak,{Bjork-2},Cui-2017,Cui-2012,Hu-jin-Zhou,Hu-jin-zhou-2} are mutually stimulated.
%
The (single-period) mean-variance formulation initiated by Markowitz \cite{Markowitz-1} is the cornerstone of modern portfolio theory and is widely used in both academic and financial industry.
The multi-period mean-variance portfolio selection is the natural extension of \cite{Markowitz-1}, which has been extensively studied.
Until 2000 and for the first time, Li-Ng \cite{Li-Duan} and Zhou-Li \cite{zhou-xunyu-2000} reported the analytical pre-commitment optimal policies for the discrete-time case and the continuous-time case, respectively.

To proceed, consider a capital market consisting of one riskless asset and $m$ risky assets
within a time horizon $N$. Let $s_k(>1)$ be a given deterministic
return of the riskless asset at time period $k$ and $e_k =
(e^1_k,\cdots,e^m_k)^T$ the vector of random returns of the $m$
risky assets at period $k$. We assume that vectors $e_k, k = 0,
1,\cdots, N-1$, are statistically independent and the only
information known about the random return vector $e_k$ is its first
two moments: its mean $\mathbb{E}(e_k) = (\mathbb{E}e^1_k,
\mathbb{E}e^2_k, \cdots , \mathbb{E}e^m_k)^T$ and its  covariance
$\mbox{Cov}(e_k)=\mathbb{E}[(e_k-\mathbb{E}e_k)(e_k-\mathbb{E}e_k)^T]$.
Clearly, $\mbox{Cov}(e_k)$ is nonnegative definite, i.e.,
$\mbox{Cov}(e_k)\succeq 0$.

Let $X_k\in \mathbb{R}$ be the wealth of the investor at the beginning
of the $k$-th period, and let $u^i_k$, $i=1,2,\cdots,m$, be the
amount invested in the $i$-th risky asset at period $k$. Then,
$X_k-\sum_{i=1}^mu_k^i$ is the amount invested in the riskless asset
at period $k$, and the wealth at the beginning of the $(k+1)$-th
period \cite{Li-Duan} is given by
\begin{eqnarray}\label{system-mean-variance}
X_{k+1}=\sum_{i=1}^ne_k^iu_k^i+\Big{(}X_k-
\sum_{i=1}^mu_k^i\Big{)}s_k=s_kX_k+O_k^T u_k,
\end{eqnarray}
where $O_k$ is the excess return vector of risky assets \cite{Li-Duan} defined as
$O_k=(O_k^1,O_k^2,\cdots,O_k^m)^T=(e_k^1-s_k,e_k^2-s_k,\cdots,e^m_k-s_k)^T$.
In this section,  we consider the case where short-selling of
stocks is allowed, i.e., $u_k^i, i=1,...,k$, can be taken values in
$\mathbb{R}$. This leads to a multi-period mean-variance
portfolio selection formulation.

Throughout this paper, we let
$\mathcal{F}_k=\sigma(e_\ell, \ell=0,1,\cdots,k-1)$, $k=0,...,N-1$.
Then, a time-inconsistent version of multi-period mean-variance problem \cite{Li-Duan} can be
formulated as follows:

\textbf{Problem (MV)}. Let $t\in \mathbb{T}$ and $x\in l^2_{\mathcal{F}}(t; \mathbb{R})$. Find $u^*\in l^2_{\mathcal{F}}(\mathbb{T}_t; \mathbb{R}^m)$ such that
\begin{eqnarray*}
{J}(t,x; u^*)=\inf_{u\in l^2_{\mathcal{F}}(\mathbb{T}_t; \mathbb{R}^m)}{J}(t,x; u).
\end{eqnarray*}
Here, $\mathbb{T}=\{0,...,N-1\}$, $\mathbb{T}_t=\{t,...,N-1\}$, and
\begin{eqnarray*}
l^2_{\mathcal{F}}(t; \mathbb{R})=\left\{\nu_t{\big{|}}
\nu_t\in \mathbb{R}\mbox{ is }\mathcal{F}_t\mbox{-measurable},~
\mathbb{E}|\nu_t|^2<\infty.
\right\},
\end{eqnarray*}
\begin{eqnarray*}
l^2_{\mathcal{F}}(\mathbb{T}_t; \mathbb{R}^{m})=\left\{\{\nu_k, k\in \mathbb{T}_t\}{\big{|}}\nu_k\in \mathbb{R}^{m}\mbox{ is }\mathcal{F}_k\mbox{-measurable}, 
\mathbb{E}|\nu_k|^2<\infty,~ k\in \mathbb{T}_t
\right\};
\end{eqnarray*}
furthermore,
\begin{eqnarray*}\label{cost-mean-variance-1}
{J}(t,x;u)=\mathbb{E}_t(X_N-\mathbb{E}_tX_N)^2-(\mu_1x+\mu_2)\mathbb{E}_tX_N,
\end{eqnarray*}
which is subject to
\begin{eqnarray*}
\left\{
\begin{array}{l}
X_{k+1}=s_kX_k+O_k^T u_k,\\[1mm]
X_t=x, ~~k\in \mathbb{T}_t
\end{array}
\right.
\end{eqnarray*}
with $\mu_1, \mu_2>0$ the trade-off parameters between the mean and
the variance of the terminal wealth.

It should be mentioned that Problem (MV) above has two unconventional features: the term $(\mu_1x)^T\mathbb{E}_tX_N$ makes $J(t,x;u)$  a state-dependent (or rank-dependent) utility, and the cost functional $J(t,x;u)$ involves the nonlinear terms of the conditional expectation of state and control variables. 
It is known now that any of the two features will ruin the time-consistency of optimal
control, namely, Bellman's principle of optimality will no longer work for Problem (MV). Note that the above model is more general than that of \cite{Ni-Zhang-Krstic2017}; in Section V of \cite{Ni-Zhang-Krstic2017}, the case without $(\mu_1x)^T\mathbb{E}_tX_N$ is dealt with.

In \cite{Li-Duan}, realizing the time-inconsistency (called nonseparability there), Li and Ng derived the optimal policy of multi-period mean-variance portfolio selection using an embedding scheme.
Note that the optimal policy of \cite{Li-Duan} is with respect to the initial pair, i.e., it is optimal only when viewed at the initial time.
This derivation is called the pre-committed optimal solution now.
By applying a pre-committed optimal control (for an initial pair), we find that it is not an optimal control for the intertemporal initial pair.
Though the pre-committed optimal solution is of some practical and theoretical values, it neglects and has not
really addressed the time-inconsistency.

In recent years, there is a surge to study the time-inconsistent optimal control together with the revisit to multi-period mean-variance portfolio selection \cite{Basak,{Bjork-2},Cui-2017,Cui-2012,Hu-jin-Zhou,Hu-jin-zhou-2,Ni-Zhang-Krstic2017,Qi-Zhang2017,Yong-1,Yong-0,Yong-2013}.
Two kinds of time-consistent equilibrium solutions are investigated in these papers, which are the
open-loop equilibrium control and the closed-loop equilibrium
strategy. To compare,
open-loop formulation is to find an open-loop equilibrium
``control", while the ``strategy" is the object of closed-loop
formulation.
Strotz's equilibrium solution \cite{Strotz} is essentially a
closed-loop equilibrium strategy, which is further elaborately
developed by Yong to the LQ optimal control \cite{Yong-1,Yong-2013}
as well as the nonlinear optimal control \cite{Yong-0,Yong-2,Wei-Yu-Yong2017}. In
contrast, open-loop equilibrium control is extensively studied by
Hu-Jin-Zhou \cite{Hu-jin-Zhou,Hu-jin-zhou-2}, Yong
\cite{Yong-2013}, Ni-Zhang-Krstic \cite{Ni-Zhang-Krstic2017}, and Qi-Zhang\cite{Qi-Zhang2017}. In particular, the closed-loop
formulation can be viewed as the extension of Bellman's dynamic
programming, and the corresponding equilibrium strategy ({if it
exists}) is derived by a backward procedure
\cite{Yong-1,Yong-2,Yong-0,Yong-2013}. Differently, the open-loop
equilibrium control is characterized via the maximum-principle-like
methodology \cite{Hu-jin-Zhou,Hu-jin-zhou-2,Ni-Zhang-Krstic2017}.

It is noted that some nondegenerate assumptions are posed in \cite{Basak,{Bjork-2},Cui-2017,Cui-2012,Hu-jin-Zhou,Hu-jin-zhou-2,Li-Duan}.
Specifically, the volatilities of the stocks in \cite{Basak,{Bjork-2},Hu-jin-Zhou,Hu-jin-zhou-2} and the return rates of the
risky securities in \cite{Cui-2017,Cui-2012,Li-Duan} are assumed to be nondegenerate, i.e., $\mbox{Cov}(e_k)\succ 0, k\in \mathbb{T}$.
%
%
To make the formulation
more practical, it is natural to consider, at least in theory, how
to generalize these results to the case where degeneracy is allowed.
In fact, mean-variance portfolio selection problems with degenerate
covariance matrices may date back to 1970s. In \cite{Buser1977} or
the ``corrected" version \cite{Ryan}, Buser \emph{et al} propose the
single-period version with possibly singular covariance matrix.
Clearly, such class of problems are more general than the classical
ones \cite{Markowitz-1}, and more consistent with the
reality.

In this note, we do not pose the nondegenerate assumption and want to find the conditions such that the time-consistent equilibrium solutions of Problem (MV) exist. This is done by using the theory developed by \cite{Ni-Li-Zhang-Krstic2018}. The rest of this paper is organized as follows. Section \ref{Sec:equilibrium-solution} gives the definitions of equilibrium solutions, whose existence is investigated in Section \ref{Section-MV}. In Section \ref{Example}, an example of \cite{Li-Duan} is revisited.

\section{Equilibrium solutions}\label{Sec:equilibrium-solution}

In the following, we shall introduce three equilibrium solutions for Problem (MV), which are the open-loop equilibrium control, feedback equilibrium strategy and mixed equilibrium solution. Note that the following notions are consistent with those of \cite{Ni-Li-Zhang-Krstic2018}. Throughout this note, Problem (MV) for the initial pair $(t,x)$ will be simply denoted as Problem (MV)$_{tx}$.

\begin{definition}
i). At stage $k\in \mathbb{T}_t$, a function $f_k(\cdot)$ is called an admissible feedback strategy (or simply a feedback strategy) if for any $\zeta\in l^2_{\mathcal{F}}(k; \mathbb{R})$, $f_k(\zeta)\in l^2_{\mathcal{F}}(k; \mathbb{R}^m)$, where
\begin{eqnarray*}
l^2_{\mathcal{F}}(k; \mathbb{R}^m)=\left\{\nu_k{\big{|}}\nu_k\in \mathbb{R}^{m}\mbox{ is }\mathcal{F}_k\mbox{-measurable},~
\mathbb{E}|\nu_k|^2<\infty.
\right\}.
\end{eqnarray*}
The set of such type of $f_k's$ is denoted by $\mathbb{F}_k$, and $\mathbb{F}_t\times \cdots\times \mathbb{F}_{N-1}$ is denoted by $\mathbb{F}_{\mathbb{T}_t}$.

ii). Let $f=(f_t,...,f_{N-1})\in \mathbb{F}_{\mathbb{T}_t}$. For $k\in \mathbb{T}_t$ and $\zeta\in l^2_{\mathcal{F}}(k; \mathbb{R})$, $f_k(\zeta)$ can be divided into two parts, namely, $f_k(\zeta)=f^c_k+f^p_k(\zeta)$; here, $f^c_k$ is the part of $f_k$ that has nothing to do with  $\zeta$ (its derivative with respect to $\zeta$ is 0), and the remainder $f_k^p(\cdot)$ is the pure-feedback-strategy part of $f_k$. Furthermore, $(f_t^p,...,f_{N-1}^p)$ is called a pure feedback strategy.

\end{definition}
%
%

\begin{definition}\label{definition-closed-loop-2}
i). A strategy $\psi\in \mathbb{F}_{\mathbb{T}_t}$ is called a feedback equilibrium strategy of Problem (MV)$_{t,x}$, if the following two points hold:

a) $\psi$ does not depend on $x$;

b) For any $k\in \mathbb{T}_t$ and $u_k\in l^2_\mathcal{F}(k; \mathbb{R}^m)$, it holds that
\begin{eqnarray}\label{defi-closed-loop}
J\big{(}k, X_k^{t,x,*}; (\psi\cdot X^{k,\psi})|_{\mathbb{T}_k}\big{)}\leq J\big{(}k, X_k^{t,x,*}; (u_k,(\psi\cdot X^{k,u_k,\psi})|_{\mathbb{T}_{k+1}})\big{)}.
\end{eqnarray}
In (\ref{defi-closed-loop}), $(\psi\cdot X^{k,\psi})|_{\mathbb{T}_k}$ and $(\psi\cdot X^{k,u_k,\psi})|_{\mathbb{T}_{k+1}}$ (with $\mathbb{T}_k=\{k,...,N-1\}, \mathbb{T}_{k+1}=\{k+1,...,N-1\}$) are given by
\begin{eqnarray*}
&&\hspace{-2.5em}(\psi\cdot X^{k,\psi})|_{\mathbb{T}_k}=\big{(}\psi_k(X^{k,\psi}_k),...,\psi_{N-1}(X^{k,\psi}_{N-1})\big{)},\\
&&\hspace{-2.5em}(\psi\cdot X^{k,u_k,\psi})|_{\mathbb{T}_{k+1}}=\big{(}\psi_{k+1}(X^{k,u_k,\psi}_{k+1}),...,\psi_{N-1}(X^{k,u_k,\psi}_{N-1})\big{)},
\end{eqnarray*}
where
$X^{k,\psi}, X^{k,u_k,\psi}$ are as follows
\begin{eqnarray}
%
%
&&\hspace{-2.25em}\label{Sec3:Defini-closed-X---2}
\left\{\begin{array}{l}
X^{k,\psi}_{\ell+1} =s_\ell X^{k,\psi}_\ell+O_\ell^T\psi_\ell(X^{k,\psi}_\ell),\\[1mm]
X^{k,\psi}_{k} = X^{t,x,*}_k,~~\ell\in \mathbb{T}_k,
\end{array}
\right.\\
&&\label{sec3:X-u---2}
\hspace{-2.25em}\left\{\begin{array}{l}
X^{k,u_k,\psi}_{\ell+1} =s_\ell X^{k,u_k,\psi}_\ell+O_\ell^T\psi_\ell(X^{k,u_k,\psi}_\ell),\\[1mm]
X^{k,u_k,\psi}_{k+1} =s_kX_k^{k,u_k,\psi}+O_\ell^Tu_k,\\[1mm]
X^{k,u_k,\psi}_{k} = X^{t,x,*}_k,~~\ell\in \mathbb{T}_{k+1}.
\end{array}\right.
\end{eqnarray}
Furthermore, in (\ref{defi-closed-loop}), (\ref{Sec3:Defini-closed-X---2}) and (\ref{sec3:X-u---2}), $X^{t,x,*}$ is given by
\begin{eqnarray*}
\left\{\begin{array}{l}
X^{t,x,*}_{k+1} =s_kX_k^{t,x,*}+O_k^T\psi_k(X^{t,x,*}_k),\\[1mm]
X^{t,x,*}_{t} = x,~~k\in \mathbb{T}_t.
\end{array}
\right.
\end{eqnarray*}

ii). Let $(\Psi, \gamma)\in l^2(\mathbb{T}_t; \mathbb{R}^{m})\times  l^2_{\mathcal{F}}(\mathbb{T}_t;\mathbb{R}^m)$ with
\begin{eqnarray*}
l^2(\mathbb{T}_t; \mathbb{R}^{m})=\left\{\{\nu_k, k\in \mathbb{T}_t\}{\Big{|}}
\nu_k\in \mathbb{R}^{m},~|\nu_k|^2<\infty, 
k\in \mathbb{T}_t
\right\}.
\end{eqnarray*}
If $\Psi$ and $\gamma$ do not depend on $x$, and $\psi$ of i) is equal to $(\Psi, \gamma)$, namely,
$$\psi_k(\xi)=\Psi_k\xi+\gamma_k, ~k\in \mathbb{T}_t,~ \xi\in l^2_{\mathcal{F}}(k;\mathbb{R}^n),$$
then $(\Psi, \gamma)$ is called a linear feedback equilibrium strategy of Problem (MV)$_{tx}$.

\end{definition}

\begin{definition}
A control $u^{t,x}\in l^2_{\mathcal{F}}(\mathbb{T}_t;\mathbb{R}^m)$ is called an open-loop equilibrium control of Problem (MV)$_{tx}$, if
\begin{eqnarray}
J\big{(}k, X_k^{t,x,*}; u^{t,x}|_{\mathbb{T}_k}\big{)}\leq J\big{(}k, X_k^{t,x,*}; (u_k,u^{t,x}|_{\mathbb{T}_{k+1}})\big{)}
\end{eqnarray}
holds for any $k\in \mathbb{T}_t$ and $u_k\in
l^2_\mathcal{F}(k; \mathbb{R}^m)$.
Here, $u^{t,x}|_{\mathbb{T}_k}$ and
$u^{t,x}|_{\mathbb{T}_{k+1}}$  are the restrictions of $u^{t,x}$
on $\mathbb{T}_k$ and $\mathbb{T}_{k+1}$, respectively; and $X^{t,x,*}$ is given by
\begin{eqnarray*}
\left\{
\begin{array}{l}
X^{t,x,*}_{k+1} = s_kX_k^{t,x,*}+O_k^Tu_k^{t,x},\\[1mm]
X^{t,x,*}_{t} = x,~~ k \in  \mathbb{T}_t.
\end{array}
\right.
\end{eqnarray*}

\end{definition}

\begin{definition}\label{definition-linear feedback}
i). A pair $(\Phi, v^{t,x}) \in l^2(\mathbb{T}_t; \mathbb{R}^{m})\times  l^2_{\mathcal{F}}(\mathbb{T}_t;\mathbb{R}^m)$ is called a mixed equilibrium solution of Problem (MV)$_{tx}$, if the following two points hold:

a) $\Phi$ does not depend on $x$, and $v^{t,x}$ depends on $x$;

b) For any $k\in \mathbb{T}_t$ and $u_k\in l^2_\mathcal{F}(k; \mathbb{R}^m)$, it holds that
\begin{eqnarray}\label{defi-linear feedback}
J\big{(}k, X_k^{t,x,*}; (\Phi\cdot X^{k,\Phi}+v^{t,x})|_{\mathbb{T}_k}\big{)}\leq J\big{(}k, X_k^{t,x,*}; (u_k,(\Phi \cdot X^{k,u_k,\Phi}+v^{t,x})|_{\mathbb{T}_{k+1}})\big{)}.
\end{eqnarray}
In (\ref{defi-linear feedback}), $(\Phi\cdot X^{k,\Phi}+v^{t,x})|_{\mathbb{T}_k}$ and $(\Phi\cdot X^{k,u_k,\Phi}+v^{t,x})|_{\mathbb{T}_{k+1}}$ are given, respectively, by
\begin{eqnarray*}
&&\hspace{-2em}\big{(}\Phi_k X_k^{k,\Phi}+v^{t,x}_k,\cdots, \Phi_{N-1} X_{N-1}^{k,\Phi}+v^{t,x}_{N-1}\big{)},\\
&&\hspace{-2em}\big{(}\Phi_{k+1} X_{k+1}^{k,u_k,\Phi}+v^{t,x}_{k+1},\cdots, \Phi_{N-1} X_{N-1}^{k,u_k,\Phi}+v^{t,x}_{N-1}\big{)},
\end{eqnarray*}
where $X^{k,\Phi}$, $X^{k,u_k,\Phi}$ are defined by
\begin{eqnarray}
%
%
&&\hspace{-2em}\label{Sec2:Defini-closed-X-1}
\left\{\begin{array}{l}
X^{k,\Phi}_{\ell+1} =\big{(}s_k+O_k^T\Phi_k\big{)}X^{k,\Phi}_k+O_k^Tv^{t,x}_k,\\[1mm]
X^{k,\Phi}_{k} = X^{t,x,*}_k,~~k\in \mathbb{T}_t,~~\ell\in \mathbb{T}_k,
\end{array}
\right.\\
&&\label{X-u}
\hspace{-2em}\left\{\begin{array}{l}
X^{k,u_k,\Phi}_{\ell+1} =\big{(}s_k+O_k^T\Phi_k \big{)}X^{k,u_k,\Phi}_k+v^{t,x}_k,\\[1mm]
X^{k,u_k,\Phi}_{k+1} =s_kX^{k,u_k,\Phi}_k+O_k^Tu_k,\\[1mm]
X^{k,u_k,\Phi}_{k} = X^{t,x,*}_k,~~\ell\in \mathbb{T}_{k+1}.
\end{array}\right.
\end{eqnarray}
The state $X^{t,x,*}_k$ in (\ref{defi-linear feedback}), (\ref{Sec2:Defini-closed-X-1}) and (\ref{X-u}) is computed via
\begin{eqnarray*}
\left\{\begin{array}{l}
X^{t,x,*}_{k+1} =\big{(}s_k+O_k^T\Phi_k\big{)}X^{t,x,*}_k+O_k^Tv^{t,x}_k,\\[1mm]
X^{t,x,*}_{t} = x,~~k\in \mathbb{T}_t.
\end{array}
\right.
\end{eqnarray*}

ii). $\Phi$ and $v^{t,x}$ in i) are called, respectively, the pure-feedback-strategy part and the open-loop-control part of the mixed equilibrium solution $(\Phi, v^{t,x})$.

iii). Letting $\Phi=0$ in i), the corresponding $v^{t,x}$ satisfying (\ref{definition-linear feedback}) is called an open-loop equilibrium control of Problem (MV)$_{tx}$.

iv). If $(\Phi, v^{t,x})$ does not depend on $x$, then it is a linear feedback equilibrium strategy of Problem (MV)$_{tx}$.

\end{definition}

\begin{remark}\label{Remark-linear feedback}
By the definition, a mixed equilibrium solution $(\Phi,v^{t,x})$ is time-consistent in the sense that  $(\Phi, v^{t,x})|_{\mathbb{T}_k}$ is a mixed  equilibrium solution for the initial pair $(k, X^{t,x,*}_k)$.
Since $(\Phi\cdot X^{k,\Phi}+v^{t,x})|_{\mathbb{T}_k}=\big{(}\Phi_k X_k^{k,\Phi}+v^{t,x}_k, (\Phi\cdot X^{k,\Phi}+v^{t,x})|_{\mathbb{T}_{k+1}}\big{)}$,
we obtain $(u_k,(\Phi\cdot X^{k,u_t,\Phi}+v^{t,x})|_{\mathbb{T}_{k+1}})\big{)}$ from $(\Phi\cdot X^{k,\Phi}+v^{t,x})|_{\mathbb{T}_k}$ by not only replacing $\Phi_k X_k^{k,\Phi}+v^{t,x}_k$ with $u_k$ but also replacing $X^{k,\Phi}$ with $X^{k,u_k, \Phi}$.
Furthermore, it is valuable to mention that $v^{t,x}$'s in both sides of (\ref{defi-linear feedback}) are the same. This is why we call $\Phi$ the pure-feedback-strategy part and $v^{t,x}$ the open-loop-control part.

\end{remark}

\section{Characterization on the equilibrium solutions}\label{Section-MV}


To solve Problem (MV)$_{tx}$, we shall transform (\ref{system-mean-variance}) into a
linear controlled system with multiplicative noises so that
the general theory \cite{Ni-Li-Zhang-Krstic2018} can work. Precisely,
define
\begin{eqnarray*}
\left\{
\begin{array}{l}
w^i_{k}=e^i_k-s_k-\mathbb{E}(e^i_k-s_k), \\
D^i_{k}=(0,\cdots,0,1,0,\cdots,0), \\
~i=1,\cdots,n,~k=0,1,\cdots,N-1,
\end{array}
\right.
\end{eqnarray*}
where the $i$-th entry of $D^i_k$ is 1.
Then, {$\{w_k=(w^1_{k},...,w^n_{k})^T, k\in \mathbb{T}\}$ is a
martingale difference sequence as $e_k, k=0,..,N-1$, are statistically independent. Furthermore,
$$\mathbb{E}_k[w_{k}w_{k}^T]=\mathbb{E}[w_{k}w_{k}^T]=\mbox{Cov}(e_{k})=(\delta^{ij}_k)_{n\times n}.$$
This leads to
\begin{eqnarray}\label{system-mean-variance-2}
\left\{
\begin{array}{l}
X_{k+1}=(s_kX_k+(\mathbb{E}O_k)^T
u_k)+\sum_{i=1}^nD^i_{k}u_kw^i_{k},\\
X_t=x,~~k\in \mathbb{T}_t.
\end{array}
\right.
\end{eqnarray}

We firstly characterize the open-loop equilibrium portfolio control of Problem (MV)$_{tx}$.

\begin{theorem}\label{Theorem-MV-open-loop}
The following statements are equivalent:

\begin{itemize}
\item[i)] For any $t\in \mathbb{T}$ and $x\in l^2_{\mathcal{F}}(t;\mathbb{R})$, Problem (MV)$_{tx}$ admits an open-loop equilibrium control.

\item[ii)] $\mathbb{E}O_k\in \mbox{Ran}(\mbox{\rm Cov}(O_k)),~ k\in \mathbb{T}_t$.
\end{itemize}

Under any of the above conditions,
\begin{eqnarray*}
v^{t,x}_k=-\widehat{\mathcal{O}}_k^\dagger \widehat{\mathcal{L}}_kX^{t,x,*}_k-\widehat{\mathcal{O}}_k^\dagger \widehat{\theta}_k,~~~k\in \mathbb{T}_t
\end{eqnarray*}
is an open-loop equilibrium control of Problem (MV)$_{tx}$, where
\begin{eqnarray*}
\left\{\begin{array}{l}
\widehat{X}^{t,x,*}_{k+1} =\big{(}s_k-O_k^T\widehat{\mathcal{O}}_k^{\dagger}\widehat{\mathcal{L}}_k \big{)}\widehat{X}^{t,x,*}_{k}-O_k^T\widehat{\mathcal{O}}^\dagger_k\widehat{\theta}_k,\\[1mm]
\widehat{X}^{t,x,*}_{t} = x,~~k\in \mathbb{T}_t,
\end{array}
\right.
\end{eqnarray*}
and
\begin{eqnarray}\label{O-L-MV}
\left\{
\begin{array}{l}
\widehat{\mathcal{O}}_{k}=\big{(}\widehat{S}_{k+1}+\widehat{T}_{k+1}\big{)}\mbox{\rm Cov}(O_k),\\[1mm] %
\widehat{\mathcal{L}}_{k}=-\frac{\mu_1}{2}s_{k+1}\cdots s_{N-1}\mathbb{E}O_k,\\[1mm]
\widehat{\theta}_{k}=-\frac{\mu_2}{2}s_{k+1}\cdots s_{N-1}\mathbb{E}O_k,\\
k\in \mathbb{T}_t
\end{array}
\right.
\end{eqnarray}
with
$$\left\{
\begin{array}{l}
\widehat{S}_k+\widehat{{T}}_k=\big{(}\widehat{S}_{k+1}+\widehat{T}_{k+1}\big{)}s_k^2-s_k\big{(}\widehat{S}_{k+1}+\widehat{T}_{k+1}\big{)}(\mathbb{E}O_k)^T
\widehat{\mathcal{O}}^\dagger_k \widehat{\mathcal{L}}_k,\\[1mm]
\widehat{S}_N+\widehat{T}_N=1, k\in \mathbb{T}_t
\end{array}
\right.
$$
\end{theorem}

\emph{Proof}. This result is proved according to Theorem 3.11 of \cite{Ni-Li-Zhang-Krstic2018}. In this case, (3.21)-(3.23) of Theorem 3.11 of \cite{Ni-Li-Zhang-Krstic2018} become to
\begin{eqnarray}\label{S-MV}
\left\{
\begin{array}{l}
\widehat{S}_{k}=s^2_k\widehat{S}_{k+1},\quad
\widehat{\mathcal{S}}_{k}=s^2_k\widehat{\mathcal{S}}_{k+1}\equiv 0, \quad\\[1mm]
\widehat{S}_N=1,~~\widehat{\mathcal{S}}_N=0,\\[1mm]
\widehat{\mathbb{O}}_k=(\mathbb{E}O_k)\widehat{\mathcal{S}}_{k+1}(\mathbb{E}O_k)^T+\sum_{i,j=1}^p\delta_k^{ij}(D_k^i)^T \widehat{S}_{k+1}D_k^j=\widehat{S}_{k+1} \mbox{Cov}(O_k)\succeq 0,\quad \\[1mm]
k\in \mathbb{T},
\end{array}
\right.
\end{eqnarray}
\begin{eqnarray}\label{T-MV}
\left\{
\begin{array}{l}
\widehat{{T}}_k=s_k^2\widehat{T}_{k+1}-s_k\big{(}\widehat{S}_{k+1}+\widehat{T}_{k+1}\big{)}(\mathbb{E}O_k)^T\widehat{\mathcal{O}}^\dagger_k \widehat{\mathcal{L}}_k,\\[1mm]
\widehat{\mathcal{T}}_k=\widehat{\mathcal{T}}_{k+1}\big{[}s_k^2-s_k(\mathbb{E}O_k)^T\widehat{\mathcal{O}}^\dagger_k \widehat{\mathcal{L}}_k\big{]}\equiv 0,\\[1mm]
\widehat{T}_N=0,~~\widehat{\mathcal{T}}_N=0,\\[1mm]
\widehat{\mathcal{O}}_k\widehat{\mathcal{O}}_k^\dagger\widehat{\mathcal{L}}_k-\widehat{\mathcal{L}}_k=0, \quad\\[1mm]
k\in \mathbb{T},
\end{array}
\right.
\end{eqnarray}
and
\begin{eqnarray}\label{pi-MV}
\left\{
\begin{array}{l}
\widehat{\pi}_k=s_k\widehat{\pi}_{k+1},\\[1mm]
\widehat{\pi}_N=-\frac{\mu_2}{2},\\[1mm]
\widehat{\mathcal{O}}_k\widehat{\mathcal{O}}_k^\dagger\widehat{\theta}_k-\widehat{\theta}_k=0, \\[1mm]
k\in \mathbb{T},
\end{array}
\right.
\end{eqnarray}
where
\begin{eqnarray*}
\left\{
\begin{array}{l}
\widehat{\mathcal{O}}_{k}=\sum_{i,j=1}^p\delta_k^{ij}({D}^i_{k})^T\big{(}\widehat{S}_{k+1}+\widehat{T}_{k+1}\big{)}{D}^j_{k}=\big{(}\widehat{S}_{k+1}+\widehat{T}_{k+1}\big{)}\mbox{Cov}(O_k),\\[1mm] %
\widehat{\mathcal{L}}_{k}=\widehat{U}_{k+1}\mathbb{E}O_k,\\[1mm]
\widehat{\theta}_{k}= \widehat{\pi}_{k+1}\mathbb{E}O_k,\\[1mm]
k\in \mathbb{T}
\end{array}
\right.
\end{eqnarray*}
with
\begin{eqnarray*}
\left\{
\begin{array}{l}
\widehat{U}_{k}=s_k\widehat{U}_{k+1},\\[1mm]
\widehat{U}_{N}=-\frac{\mu_1}{2},\\[1mm]
k\in \mathbb{T}.
\end{array}
\right.
\end{eqnarray*}
From (\ref{S-MV}), (\ref{T-MV}) and (\ref{pi-MV}), we have
\begin{eqnarray*}\label{}
\left\{
\begin{array}{l}
\widehat{S}_k+\widehat{{T}}_k=\big{(}\widehat{S}_{k+1}+\widehat{T}_{k+1}\big{)}s_k^2-s_k\big{(}\widehat{S}_{k+1}+\widehat{T}_{k+1}\big{)}(\mathbb{E}O_k)^T\widehat{\mathcal{O}}^\dagger_k \widehat{\mathcal{L}}_k,\\[1mm]
%
%
\widehat{S}_N+\widehat{T}_N=1,\\[1mm]
\widehat{\mathcal{O}}_k\widehat{\mathcal{O}}_k^\dagger\widehat{\mathcal{L}}_k-\widehat{\mathcal{L}}_k=0, \\[1mm]
k\in \mathbb{T}.
\end{array}
\right.
\end{eqnarray*}
By some calculations, we have
\begin{eqnarray*}
&&\hspace{-3em}\widehat{S}_k+\widehat{{T}}_k=\big{(}\widehat{S}_{k+1}+\widehat{T}_{k+1}\big{)}s_k^2+s_ks_{k+1}\cdots s_{N-1}\frac{\mu_1}{2}\big{(}\widehat{S}_{k+1}+\widehat{T}_{k+1}\big{)}\big{(}\widehat{S}_{k+1}+\widehat{T}_{k+1}\big{)}^\dagger (\mathbb{E}O_k)^T\big{[}\mbox{Cov}(O_k)\big{]}^\dagger \mathbb{E}O_k,
\end{eqnarray*}
where $
\big{(}\widehat{S}_{k+1}+\widehat{T}_{k+1}\big{)}^\dagger$equals to $\big{(}\widehat{S}_{k+1}+\widehat{T}_{k+1}\big{)}^{-1}$ if $\widehat{S}_{k+1}+\widehat{T}_{k+1}\neq 0$; otherwise, it will be $0$.
Therefore,
\begin{eqnarray*}
\big{(}\widehat{S}_{k+1}+\widehat{T}_{k+1}\big{)}\big{(}\widehat{S}_{k+1}+\widehat{T}_{k+1}\big{)}^\dagger=\left\{
\begin{array}{ll}
1,&\widehat{S}_{k+1}+\widehat{T}_{k+1}\neq 0,\\
0,&\widehat{S}_{k+1}+\widehat{T}_{k+1}=0.
\end{array}
\right.
\end{eqnarray*}
As $s_k>1, \mu_1>0,\mbox{Cov}(O_k)\succeq 0, k\in \mathbb{T}$, and $\widehat{S}_N+\widehat{T}_N=1$, it follows that
$\widehat{S}_k+\widehat{T}_k>0,~~~k\in \mathbb{T}$.
This together with $\widehat{\pi}_k\neq 0, \widehat{U}_k\neq 0, k\in \mathbb{T}$ implies that the solvability of (\ref{S-MV})-(\ref{pi-MV}) is
equivalent to the property $\mathbb{E}O_k\in \mbox{Ran}(\mbox{Cov}(O_k)), k\in \mathbb{T}$. Then, the proof is completed by following Theorem 3.11 of \cite{Ni-Li-Zhang-Krstic2018}.  \hfill $\square$



\begin{corollary}
Let $\mbox{\rm Cov}(O_k)\succ0, k\in \mathbb{T}$. Then, for any $t\in \mathbb{T}$ and $x\in l^2_{\mathcal{F}}(t; \mathbb{R})$, Problem (MV) admits a unique open-loop equilibrium control, which is given by
$v^{t,x}_k=-\widehat{\mathcal{O}}_k^{-1} \widehat{\mathcal{L}}_kX^{t,x,*}_k-\widehat{\mathcal{O}}_k^{-1} \widehat{\theta}_k$,
$k\in \mathbb{T}_t$
with
\begin{eqnarray*}
\left\{\begin{array}{l}
\widehat{X}^{t,x,*}_{k+1} =\big{(}s_k-O_k^T\widehat{\mathcal{O}}_k^{-1}\widehat{\mathcal{L}}_k \big{)}\widehat{X}^{t,x,*}_{k}-O_k^T\widehat{\mathcal{O}}^{-1}_k\widehat{\theta}_k,\\[1mm]
\widehat{X}^{t,x,*}_{t} = x,~~k\in \mathbb{T}_t,
\end{array}
\right.
\end{eqnarray*}
and $\widehat{O}_k, \widehat{L}_k, \widehat{\theta}_k, k\in \mathbb{T}_t$, are given in (\ref{O-L-MV}).

\end{corollary}

\emph{Proof}. This follows from Theorem 3.9 of \cite{Ni-Li-Zhang-Krstic2018} and Theorem \ref{Theorem-MV-open-loop}. \hfill $\square$

The following result is based on Theorem 3.13 of \cite{Ni-Li-Zhang-Krstic2018}, which is about the feedback equilibrium strategy of Problem (MV)$_{tx}$.

\begin{theorem}\label{Theorem-MV-feedback}
The following statements are equivalent:

\begin{itemize}
\item[i)] For any $t\in \mathbb{T}$ and $x\in l_\mathcal{F}^2(t;\mathbb{R})$, Problem (MV)$_{tx}$ admits a feedback equilibrium strategy.

\item[ii)] The difference equations
\begin{eqnarray}\label{S-MV-feedback-1}
\left\{
\begin{array}{l}
\widetilde{S}_k=\widetilde{S}_{k+1}\Big{[}s_k^2-2s_k(\mathbb{E}O_k)^T\widetilde{\mathbb{O}}_k^\dagger \widetilde{\mathbb{L}}_k+\widetilde{\mathbb{L}}^T_k\widetilde{\mathbb{O}}_k^\dagger \mathbb{E}(O_kO_k^T)\widetilde{\mathbb{O}}_k^\dagger \widetilde{\mathbb{L}}_k\Big{]},\\[1mm]
\widetilde{\mathcal{S}}_k=s_k^2\widetilde{\mathcal{S}}_{k+1}-2s_k\widetilde{\mathcal{S}}_{k+1}(\mathbb{E}O_k)^T\widetilde{\mathbb{O}}_k^\dagger \widetilde{\mathbb{L}}_k+\widetilde{\mathbb{L}}^T_k\widetilde{\mathbb{O}}_k^\dagger \widetilde{\mathbb{L}}_k,\\[1mm]
\widetilde{S}_{N}=1,~~\widetilde{\mathcal{S}}_N=0,\quad\\[1mm]
\widetilde{\mathbb{O}}_k\succeq 0,~~
\widetilde{\mathbb{O}}_k\widetilde{\mathbb{O}}_k^\dagger\widetilde{\mathbb{L}}_k=\widetilde{\mathbb{L}}_k,\\[1mm]
k\in \mathbb{T},
\end{array}\right.
\end{eqnarray}
and
\begin{eqnarray}\label{pi-MV-feedback}
\hspace{-0.5em}\left\{
\begin{array}{l}
\widetilde{\pi}_{k}=-\widetilde{\beta}_{k}\widetilde{{\mathbb{O}}}^\dagger_k\widetilde{\theta}_k+\big{(}s_{k}-(\mathbb{E}O_k)^T\widetilde{\mathbb{O}}_k^\dagger\widetilde{\mathbb{L}}_k\big{)}^T \widetilde{\pi}_{k+1},\\[1mm]
\pi_{N}=-\frac{\mu_2}{2},\\[1mm]
\widetilde{\mathbb{O}}_k\widetilde{\mathbb{O}}_k^\dagger\widetilde{\theta}_k=\widetilde{\theta}_k,\\[1mm]
k\in \mathbb{T}
\end{array}
\right.
\end{eqnarray}
are solvable, namely,
\begin{eqnarray*}
\left\{
\begin{array}{l}
\widetilde{\mathbb{O}}_k\succeq 0,~~\\[1mm]
\widetilde{\mathbb{O}}_k\widetilde{\mathbb{O}}_k^\dagger\widetilde{\mathbb{L}}_k=\widetilde{\mathbb{L}}_k,\\[1mm]
\widetilde{\mathbb{O}}_k\widetilde{\mathbb{O}}_k^\dagger\widetilde{\theta}_k=\widetilde{\theta}_k,~~\\[1mm]
k\in \mathbb{T},
\end{array}
\right.
\end{eqnarray*}
holds, where
\begin{eqnarray}\label{O-pi-theta-Mv}
\left\{
\begin{array}{l}
\widetilde{\mathbb{O}}_{k}=\widetilde{\mathcal{S}}_{k+1}\mathbb{E}O_k(\mathbb{E}O_k)^T+\widetilde{S}_{k+1}\mbox{Cov}(O_k),\\[1mm]
\widetilde{\mathbb{L}}_{k}=\big{(}s_k\widetilde{\mathcal{S}}_{k+1}+\widetilde{U}_{k+1}\big{)}\mathbb{E}O_k,\\[1mm]
\widetilde{\theta}_{k}=\widetilde{\pi}_{k+1}\mathbb{E}O_k,\\[1mm]
k\in \mathbb{T}
\end{array}
\right.
\end{eqnarray}
with
\begin{eqnarray*}
\left\{
\begin{array}{l}
\widetilde{U}_{k}=\big{(}s_k-(\mathbb{E}O_k)^T\widetilde{\mathbb{O}}_k^\dagger\widetilde{\mathbb{L}}_k\big{)}\widetilde{U}_{k+1},\\[1mm]
\widetilde{U}_{N}=-\frac{\mu_1}{2},\\[1mm] k\in \mathbb{T},
\end{array}
\right.
\end{eqnarray*}
and
\begin{eqnarray*}
\left\{
\begin{array}{l}
\widetilde{\beta}_{k}=s_k\widetilde{\mathcal{S}}_{k+1}(\mathbb{E}O_k)^T-\widetilde{\mathbb{L}}_k^T\widetilde{\mathbb{O}}_k^\dagger \big{[}\widetilde{\mathcal{S}}_{k+1}\mathbb{E}O_k(\mathbb{E}O_k)^T\\[1mm]
\hphantom{\widetilde{\beta}_{k}=}+\widetilde{S}_{k+1}\mbox{Cov}(O_k)\big{]},\\[1mm]
k\in \mathbb{T}.
\end{array}
\right.
\end{eqnarray*}
\end{itemize}
\end{theorem}

\begin{theorem}
If $\mathbb{E}O_k\in \mbox{Ran}(O_k), k\in \mathbb{T}$, then (\ref{S-MV-feedback-1}) and (\ref{pi-MV-feedback}) are solvable, and for any $t\in \mathbb{T}$ and $x\in l_\mathcal{F}^2(t;\mathbb{R})$, Problem (MV)$_{tx}$ admits a feedback equilibrium strategy
$(\Phi^t, v^t)$ with
\begin{eqnarray*}
%
\Phi^t=\{-\widetilde{\mathbb{O}}^{\dagger}_{k} \widetilde{\mathbb{L}}_{k},~k\in \mathbb{T}_t\},~~~
v^t=\{-\widetilde{\mathbb{O}}^{\dagger}_{k} \widetilde{\theta}_{k},~k\in \mathbb{T}_t\},
\end{eqnarray*}
where $\widetilde{\mathbb{O}}_{k}, \widetilde{\mathbb{L}}_{k}, \widetilde{\theta}_{k}, k\in \mathbb{T}_t$, are given in (\ref{O-pi-theta-Mv}).

\end{theorem}

\emph{Proof}.  Note that (\ref{S-MV-feedback-1}) can be equivalently rewritten as
\begin{eqnarray}\label{S-MV-feedback-2}
\left\{
\begin{array}{l}
\widetilde{S}_k=\widetilde{{S}}_{k+1}\Big{[}\big{(}s_k-(\mathbb{E}O_k)^T\widetilde{\mathbb{O}}_k^\dagger\widetilde{\mathbb{L}}_k\big{)}^2+\widetilde{\mathbb{L}}^T_k\widetilde{\mathbb{O}}_k^\dagger \mbox{Cov}(O_k)\widetilde{\mathbb{O}}_k^\dagger \widetilde{\mathbb{L}}_k\Big{]},\\[1mm]
\widetilde{\mathcal{S}}_k=\widetilde{\mathcal{S}}_{k+1}\big{(}s_k-(\mathbb{E}O_k)^T\widetilde{\mathbb{O}}_k^\dagger\widetilde{\mathbb{L}}_k\big{)}^2 +\widetilde{S}_{k+1}\widetilde{\mathbb{L}}^T_k\widetilde{\mathbb{O}}_k^\dagger \mbox{Cov}(O_k)\widetilde{\mathbb{O}}_k^\dagger \widetilde{\mathbb{L}}_k,\\[1mm]
\widetilde{S}_{N}=1,~~\widetilde{\mathcal{S}}_N=0,\quad
~~\\[1mm]
\widetilde{\mathbb{O}}_k\succeq 0, \widetilde{\mathbb{O}}_k\widetilde{\mathbb{O}}_k^\dagger\widetilde{\mathbb{L}}_k=\widetilde{\mathbb{L}}_k,\\[1mm]
k\in \mathbb{T}.
\end{array}\right.
\end{eqnarray}
Clearly,  $\widetilde{S}_k\geq \widetilde{\mathcal{S}}_k\geq 0, k\in \mathbb{T}$.

Let $\mathbb{E}O_k\in \mbox{Ran}(\mbox{Cov}(O_k)), k\in \mathbb{T}$. We now prove that (\ref{S-MV-feedback-1}) and (\ref{pi-MV-feedback}) are solvable. For a generic $k\in \mathbb{T}$, we prove the conclusion by the following two cases.

Case 1: $\widetilde{S}_{k}>0$. This implies
$\widetilde{S}_\ell>0$,  $\big{(}s_\ell-(\mathbb{E}O_\ell)^T\widetilde{\mathbb{O}}_\ell^\dagger\widetilde{\mathbb{L}}_\ell\big{)}^2+\widetilde{\mathbb{L}}^T_\ell\widetilde{\mathbb{O}}_\ell^\dagger \mbox{Cov}(O_\ell)\widetilde{\mathbb{O}}_\ell^\dagger \widetilde{\mathbb{L}}_\ell>0$,
$\ell\in \mathbb{T}_k$.
If $s_k\widetilde{\mathcal{S}}_{k+1}+\widetilde{U}_{k+1}=0$, then $\widetilde{\mathbb{L}}_k=0$ and (\ref{S-MV-feedback-2}) is solvable at $k$. To the end of this paragraph, we assume $s_k\widetilde{\mathcal{S}}_{k+1}+\widetilde{U}_{k+1}\neq 0$. If further $\widetilde{\mathcal{S}}_{k+1}=0$, then $\widetilde{\mathbb{L}}_k\in \mbox{Ran}(\widetilde{\mathbb{O}}_k)$ and (\ref{S-MV-feedback-2}) is solvable at $k$. On the other hand, let $\widetilde{\mathcal{S}}_{k+1}\neq 0$. As $\mathbb{E}O_k\in \mbox{Ran}(\mbox{Cov}(O_k))$, there exists $\xi\in \mathbb{R}^n$ such that $\mbox{Cov}(O_k)\xi=\mathbb{E}O_k$. Furthermore, $\xi^T\mbox{Cov}(O_k)\xi=\xi^T\mathbb{E}O_k\geq 0$. Then,
\begin{eqnarray*}
&&\hspace{-2em}\widetilde{\mathbb{O}}_k\frac{s_k\widetilde{\mathcal{S}}_{k+1}+\widetilde{U}_{k+1}}{\widetilde{\mathcal{S}}_{k+1}\xi^T\mathbb{E}O_k+\widetilde{S}_{k+1}} \xi =\frac{s_k\widetilde{\mathcal{S}}_{k+1}+\widetilde{U}_{k+1}}{\widetilde{\mathcal{S}}_{k+1}\xi^T\mathbb{E}O_k+ \widetilde{S}_{k+1}}\big{(}\widetilde{\mathcal{S}}_{k+1}\mathbb{E}O_k(\mathbb{E}O_k)^T\xi +\widetilde{S}_{k+1}\mbox{Cov}(O_k)\xi\big{)}\\
&&\hspace{-2em}\hphantom{\widetilde{\mathbb{O}}_k\frac{s_k\widetilde{\mathcal{S}}_{k+1}+\widetilde{U}_{k+1}}{\widetilde{\mathcal{S}}_{k+1}\xi^T\mathbb{E}O_k+\widetilde{S}_{k+1}} \xi }=\big{(}s_k\widetilde{\mathcal{S}}_{k+1}+\widetilde{U}_{k+1}\big{)}\mathbb{E}O_k
=\widetilde{\mathbb{L}}_k.
\end{eqnarray*}
Hence, $\widetilde{\mathbb{O}}_k\widetilde{\mathbb{O}}_k^\dagger\widetilde{\mathbb{L}}_k=\widetilde{\mathbb{L}}_k$, and (\ref{S-MV-feedback-2}) is solvable at $k$. By a similar procedure, we can prove the solvability of (\ref{pi-MV-feedback}) at $k$.

Case 2: $\widetilde{S}_{k}=0$. If $\widetilde{S}_{k+1}>0$, the proof of the solvability of (\ref{pi-MV-feedback}) and (\ref{S-MV-feedback-2}) at $k$ is similar to that of Case 1, and is omitted here. If $\widetilde{S}_{k+1}=0, \widetilde{S}_{k+2}>0$, we have $\widetilde{\mathbb{O}}_k=0$ and
\begin{eqnarray}\label{S-MV-feedback-4}
&&\hspace{-3em}s_{k+1}-(\mathbb{E}O_{k+1})^T\widetilde{\mathbb{O}}_{k+1}^\dagger\widetilde{\mathbb{L}}_{k+1}=\widetilde{\mathbb{L}}^T_{k+1}\widetilde{\mathbb{O}}_{k+1}^\dagger \mbox{Cov}(O_{k+1})\widetilde{\mathbb{O}}_{k+1}^\dagger \widetilde{\mathbb{L}}_{k+1}=0,
\end{eqnarray}
which further implies $\widetilde{U}_{k+1}=0$ and $\widetilde{\mathbb{L}}_k=0$; hence, (\ref{S-MV-feedback-2}) is solvable at $k$. Furthermore, (\ref{S-MV-feedback-4}) implies
$\widetilde{\mathbb{L}}^T_{k+1}\widetilde{\mathbb{O}}_{k+1}^\dagger \mbox{Cov}(O_{k+1})=0$.
From this and (\ref{S-MV-feedback-4}) we have $\widetilde{\beta}_{k+1}=0$ and $\widetilde{\pi}_{k+1}=0$, and hence (\ref{pi-MV-feedback}) is solvable at $k$ under the condition of $\widetilde{S}_{k+1}=0, \widetilde{S}_{k+2}>0$. Finally, if $\widetilde{S}_{k+1}=0, \widetilde{S}_{k+2}=0$, we must have some $\tau>k+1$ such that $\widetilde{S}_\tau=0, \widetilde{S}_{\tau+1}>0$. Similar to the comments below (\ref{S-MV-feedback-4}), we have $\widetilde{U}_\tau=\widetilde{\beta}_\tau=\widetilde{\pi}_\tau=0$, which implies $\widetilde{U}_{k+1}=0$ and the solvability of (\ref{S-MV-feedback-2}) at $k$. As $\widetilde{\mathbb{O}}_\ell=0, \ell\in \{k,...,\tau-1 \}$, it follows that
$\widetilde{\pi}_{k+1}=\big{(}s_{k+1}-(\mathbb{E}O_{k+1})^T\widetilde{\mathbb{O}}_{k+1}^\dagger\widetilde{\mathbb{L}}_{k+1}\big{)}^T\cdots \big{(}s_{\tau-1}-(\mathbb{E}O_{\tau-1})^T\widetilde{\mathbb{O}}_{\tau-1}^\dagger\widetilde{\mathbb{L}}_{\tau-1}\big{)}^T\widetilde{\pi}_\tau=0$.
Hence, (\ref{pi-MV-feedback}) is solvable at $k$.

In summary, for a generic $k\in \mathbb{T}$, we have proved the solvability of (\ref{pi-MV-feedback}) and (\ref{S-MV-feedback-2}) at $k$, namely, $\mathbb{}\widetilde{\mathbb{O}}_k\widetilde{\mathbb{O}}_k^\dagger \widetilde{\mathbb{L}}_k=\widetilde{\mathbb{L}}_k$ and $\mathbb{}\widetilde{\mathbb{O}}_k\widetilde{\mathbb{O}}_k^\dagger \widetilde{\theta}_k=\widetilde{\theta}_k$. From Theorem 3.13 of \cite{Ni-Li-Zhang-Krstic2018} and Theorem \ref{Theorem-MV-feedback}, we can complete the proof.  \hfill $\square$

\begin{theorem}
Let $\mbox{Cov}(O_k)\succ0, k\in \mathbb{T}$. Then, for any $t\in \mathbb{T}$ and $x\in l^2_{\mathcal{F}}(t;\mathbb{R})$, Problem (MV)$_{tx}$ admits a unique feedback equilibrium strategy, which is given by  $(\Phi^t, v^t)$ with
$\Phi^t=\{-\widetilde{\mathbb{O}}^{-1}_{k} \widetilde{\mathbb{L}}_{k},~k\in \mathbb{T}_t\}$,
$v^t=\{-\widetilde{\mathbb{O}}^{-1}_{k} \widetilde{\theta}_{k},~k\in \mathbb{T}_t\}$,
where $\widetilde{\mathbb{O}}_{k}, \widetilde{\mathbb{L}}_{k}, \widetilde{\theta}_{k}, k\in \mathbb{T}_t$, are given in (\ref{O-pi-theta-Mv}).
\end{theorem}

\emph{Proof.} In this case, (\ref{S-MV-feedback-1}) and (\ref{pi-MV-feedback}) are solvable.  From (\ref{S-MV-feedback-2}), we know that $\widetilde{S}_{k}>0, k\in \mathbb{T}$. In fact, suppose $\widetilde{S}_{k_0}=0$ and $\widetilde{S}_{k_0+1}\neq 0$ for some $k_0\in \mathbb{T}$, then
\begin{eqnarray*}
s_{k_0}-(\mathbb{E}O_{k_0})^T\widetilde{\mathbb{O}}_{k_0}^\dagger\widetilde{\mathbb{L}}_{k_0}=\widetilde{\mathbb{L}}^T_{k_0}\widetilde{\mathbb{O}}_{k_0}^\dagger \mbox{Cov}(O_{k_0})\widetilde{\mathbb{O}}_{k_0}^\dagger \widetilde{\mathbb{L}}_{k_0}=0.
\end{eqnarray*}
As $\mbox{Cov}(O_{k_0})\succ0$, it follows $\widetilde{\mathbb{O}}_{k_0}^\dagger \widetilde{\mathbb{L}}_{k_0}=0$, which implies $0=s_{k_0}-(\mathbb{E}O_{k_0})^T\widetilde{\mathbb{O}}_{k_0}^\dagger\widetilde{\mathbb{L}}_{k_0}=s_{k_0}$. This is impossible, and thus $\widetilde{S}_{k}>0, k\in \mathbb{T}$. 
%
%
%
%
%
%
%
%
Furthermore, we have $\widetilde{\mathbb{O}}_k\succ 0, k\in \mathbb{T}$. This completes the proof. \hfill $\square$

We now consider the mixed equilibrium portfolio solution. In this case,  (3.26)-(3.28) of \cite{Ni-Li-Zhang-Krstic2018} read as
\begin{eqnarray}\label{S-MV-mixed}
\left\{
\begin{array}{l}
{S}_{k}={S}_{k+1}\Big{[}\big{(}s_k+(\mathbb{E}O_k)^T\Phi_k\big{)}^2+\Phi_k^T\mbox{Cov}(O_k)\Phi_k\Big{]},\\[1mm]
{\mathcal{S}}_{k}={\mathcal{S}}_{k+1}\big{(}s_k+(\mathbb{E}O_k)^T\Phi_k\big{)}^2+S_{k+1}\Phi_k^T\mbox{Cov}(O_k)\Phi_k,  \\[1mm]
{S}_N=1,~~{\mathcal{S}}_N=0,\\
{\mathbb{O}}_{k}={\mathcal{S}}_{k+1}\mathbb{E}O_k(\mathbb{E}O_k)^T+{S}_{k+1}\mbox{Cov}(O_k)\succeq 0,\\[1mm]
k\in \mathbb{T},
\end{array}
\right.
\end{eqnarray}
\begin{eqnarray}\label{T-MV-mixed}
\left\{
\begin{array}{l}
T_k=S_{k+1}\Big{[}\big{(}s_k+(\mathbb{E}O_k)^T\Phi_k\big{)}(\mathbb{E}O_k)^T+\Phi_k^T\mbox{Cov}(O_k)\Big{]}\big{(}-\mathcal{O}_k^\dagger \mathcal{L}_k-\Phi_k \big{)}\\[1mm]
\hphantom{T_k=}+T_{k+1}\Big{[}\big{(}s_k+(\mathbb{E}O_k)^T\Phi_k\big{)}\big{(}s_k-(\mathbb{E}O_k)^T\mathcal{O}_k^\dagger \mathcal{L}_k \big{)}-\Phi_k^T\mbox{Cov}(O_k)\mathcal{O}_k^\dagger \mathcal{L}_k \Big{]},\\[1mm]
\mathcal{T}_k=\Big{[}\mathcal{S}_{k+1}\big{(}s_k+(\mathbb{E}O_k)^T\Phi_k\big{)}(\mathbb{E}O_k)^T+{S}_{k+1}\Phi_k^T\mbox{Cov}(O_k)\Big{]}\big{(}-\mathcal{O}_k^\dagger \mathcal{L}_k -\Phi_k\big{)}\\[1mm]
\hphantom{T_k=}+\Big{[}\mathcal{T}_{k+1}\big{(}s_k+(\mathbb{E}O_k)^T\Phi_k\big{)}\big{(}s_k-(\mathbb{E}O_k)^T\mathcal{O}_k^\dagger \mathcal{L}_k \big{)}-T_{k+1}\Phi_k^T\mbox{Cov}(O_k)\mathcal{O}_k^\dagger \mathcal{L}_k \Big{]},\\[1mm]
T_{N}=0,~~\mathcal{T}_{N}=0,\\[1mm]
\mathcal{O}_k \mathcal{O}_k^\dagger\mathcal{L}_k=\mathcal{L}_k,\\[1mm]
k\in \mathbb{T},
\end{array}
\right.
\end{eqnarray}
and
\begin{eqnarray}\label{pi-MV-mixed}
\left\{
\begin{array}{l}
{\pi}_{k}=-{\beta}_{k}{{\mathcal{O}}}^\dagger_k{\theta}_k+\big{(}s_{k}+(\mathbb{E}O_k)^T\Phi_k\big{)}^T {\pi}_{k+1}, \\[1mm]
\pi_{N}=-\frac{\mu_2}{2},\\[1mm]
{\mathcal{O}}_k{\mathcal{O}}_k^\dagger{\theta}_k={\theta}_k,\\[1mm]
k\in \mathbb{T},
\end{array}
\right.
\end{eqnarray}
where
\begin{eqnarray*}
\left\{
\begin{array}{l}
\mathcal{O}_{k}=\big{(}\mathcal{S}_{k+1}+\mathcal{T}_{k+1}\big{)}\mathbb{E}O_k(\mathbb{E}O_k)^T+\big{(}{S}_{k,k+1}+T_{k,k+1}\big{)}\mbox{Cov}(O_k),\\[1mm]
\mathcal{L}_{k}=s_k\big{(}\mathcal{S}_{k+1}+\mathcal{T}_{k+1}\big{)}\mathbb{E}O_k+U_{k+1}\mathbb{E}O_k,\\[1mm]
\theta_{k}=\pi_{k+1}\mathbb{E}O_k,\\[1mm]
k\in \mathbb{T},
\end{array}
\right.
\end{eqnarray*}
with
${U}_{k}=\big{(}s_k+(\mathbb{E}O_k)^T\Phi_k\big{)}{U}_{k+1}$,
${U}_{N}=-\frac{\mu_1}{2}$,
$k\in \mathbb{T}$,
and
${\beta}_{k}=\big{(}{\mathcal{S}}_{k+1}+\mathcal{T}_{k+1}\big{)}\big{(}s_k+\Phi_k^T\mathbb{E}O_k\big{)}(\mathbb{E}O_k)^T +\big{(}S_{k+1}+T_{k+1}\big{)}\Phi_k^T\mbox{Cov}(O_k)$, $k\in \mathbb{T}$.
From (\ref{S-MV-mixed})-(\ref{T-MV-mixed}), we obtain
\begin{eqnarray*}
\left\{
\begin{array}{l}
S_{k}+T_k=(S_{k+1}+T_{k+1})\Big{[}\big{(}s_k+(\mathbb{E}O_k)^T\Phi_k\big{)}\big{(}s_k-(\mathbb{E}O_k)^T\mathcal{O}_k^\dagger \mathcal{L}_k\big{)}-\Phi_k^T\mbox{Cov}(O_k)\mathcal{O}_k^\dagger \mathcal{L}_k\Big{]},\\[1mm]
\mathcal{S}_{k}+\mathcal{T}_k=(\mathcal{S}_{k+1}+\mathcal{T}_{k+1})\big{(}s_k+(\mathbb{E}O_k)^T\Phi_k\big{)}\big{(}s_k-(\mathbb{E}O_k)^T\mathcal{O}_k^\dagger \mathcal{L}_k\big{)}\\[1mm]
\hphantom{\mathcal{S}_{k}+\mathcal{T}_k=}-(S_{k+1}+T_{k+1})\Phi_k^T\mbox{Cov}(O_k)\mathcal{O}_k^\dagger \mathcal{L}_k,\\[1mm]
k\in \mathbb{T}.
\end{array}
\right.
\end{eqnarray*}

Note that for any $k\in \mathbb{T}$, $\mathbb{O}_k\succeq 0$, (\ref{S-MV-mixed}) is solvable. By Theorem 3.14 of \cite{Ni-Li-Zhang-Krstic2018}, the following result is straightforward.

\begin{proposition}\label{Proposition-MV-mixed}
Assume that there exists $\Phi\in l^2(\mathbb{T}; \mathbb{R}^{m})$ such that (\ref{T-MV-mixed}) (\ref{pi-MV-mixed}) are solvable. For any $t\in \mathbb{T}$ and $x\in l^2_\mathcal{F}(t; \mathbb{R})$, let
$$v^{t,x}_k=-({\mathcal{O}}_k^{\dagger}{\mathcal{L}}_k+\Phi_k)X^{t,x,*}_k-\mathcal{O}^\dagger_k\theta_k, ~~k\in \mathbb{T}_t,$$
where
$$
\left\{
\begin{array}{l}
{X}^{t,x,*}_{k+1} =\big{(}s_k-O_k^T{\mathcal{O}}_k^{\dagger}{\mathcal{L}}_k \big{)}{X}^{t,x,*}_{k}-O_k^T{\mathcal{O}}^\dagger_k{\theta}_k,\\[1mm]
{X}^{t,x,*}_{t} = x,~~k\in \mathbb{T}_t
\end{array}
\right.
$$
Then, $(\Phi|_{\mathbb{T}_t}, v^{t,x})$ is a mixed equilibrium solution of Problem (MV)$_{tx}$.

\end{proposition}

\section{Revisiting an example of \cite{Li-Duan}}\label{Example}

Consider a multi-period mean-variance portfolio selection problem (Example 2 of \cite{Li-Duan}). A capital market consists of one riskless asset and three risky assets
over  a finite time horizon $N=4$, and the parameters of the model are as follows
\begin{eqnarray*}
&&x=1,~~s_k=1.04,~~ \mathbb{E}e_k^1=1.162,~~\mathbb{E}e_k^2=1.246,~~\\
&&\mathbb{E}e_k^3=1.228,~~k=0,1,2,3,
\end{eqnarray*}
and
the covariance of $e_k=(e_k^1, e_k^2,e^3_k)^T$ is
\begin{eqnarray*}
\mbox{Cov}(e_k)=\left[
\begin{array}{ccc}
0.0146&0.0187&0.0145\\[0.5mm]
0.0187&0.0854&0.0104\\[0.5mm]
0.0145&0.0104&0.0289
\end{array}
\right]\succ0,~k=0,1,2,3.
\end{eqnarray*}
In this paper, we assume $\mu_1=\mu_2=1$. Clearly,
$$\mathbb{E}O_k=(0.1220, ~0.2060, ~0.1880)^T, k=0,1,2,3.$$

As $\mbox{Cov}(O_k)=\mbox{Cov}(e_k)\succ 0, k=0,1,2,3$, this Problem (MV) will have a unique open-loop equilibrium control and a unique feedback equilibrium strategy.
In what follows, we will compute the equilibrium solutions, respectively.

\textbf{\underline{Open-loop equilibrium control}}

By (\ref{O-L-MV}) and some calculations, we have
\begin{eqnarray*}
&&\hspace{-2em}-\widehat{\mathcal{O}}^\dagger_3\widehat{\mathcal{L}}_3=\left[
\begin{array}{c}
    0.4739\\
    0.7689\\   2.7381
\end{array}
\right],~
-\widehat{\mathcal{O}}^\dagger_2\widehat{\mathcal{L}}_2=\left[
\begin{array}{c}
    0.2676\\
    0.4341\\
    1.5461
\end{array}
\right],~\\
&&\hspace{-2em}-\widehat{\mathcal{O}}^\dagger_1\widehat{\mathcal{L}}_1=\left[
\begin{array}{c}
    0.1842\\
    0.2988\\
    1.0643
\end{array}
\right],~
-\widehat{\mathcal{O}}^\dagger_0\widehat{\mathcal{L}}_0=\left[
\begin{array}{c}
    0.1391\\
    0.2257\\
    0.8038
\end{array}
\right],\\
&&\hspace{-2em}-\widehat{\mathcal{O}}^\dagger_3\widehat{\theta}_3=\left[
\begin{array}{c}
    0.4739\\
    0.7689\\
    2.7381
\end{array}
\right],~
-\widehat{\mathcal{O}}^\dagger_2\widehat{\theta}_2=\left[
\begin{array}{c}
    0.2676\\
    0.4341\\
    1.5461
\end{array}
\right],~\\
&&\hspace{-2em}-\widehat{\mathcal{O}}^\dagger_1\widehat{\theta}_1=\left[
\begin{array}{c}
    0.1842\\
    0.2988\\
    1.0643
\end{array}
\right],~
-\widehat{\mathcal{O}}^\dagger_0\widehat{\theta}_0=\left[
\begin{array}{c}
    0.1391\\
    0.2257\\
    0.8038
\end{array}
\right],
\end{eqnarray*}
Then, the unique open-loop equilibrium portfolio control for the initial pair $(0,x)$ is given by
$$v^{0,x}_k=-\widehat{\mathcal{O}}_k^\dagger \widehat{\mathcal{L}}_kX^{0,x,*}_k-\widehat{\mathcal{O}}_k^\dagger \widehat{\theta}_k$$
with
$$
\left\{
\begin{array}{l}
\widehat{X}^{0,x,*}_{k+1} =\big{(}s_k-O_k^T\widehat{\mathcal{O}}_k^{\dagger}\widehat{\mathcal{L}}_k \big{)}\widehat{X}^{0,x,*}_{k}-O_k^T\widehat{\mathcal{O}}^\dagger_k\widehat{\theta}_k,\\
\widehat{X}^{0,x,*}_{0} = x, k=0,1,2,3.
\end{array}
\right.
$$

\textbf{\underline{Feedback equilibrium strategy}}

By (\ref{S-MV-feedback-1}), (\ref{pi-MV-feedback}) and some calculations, we have
\begin{eqnarray*}
&&\hspace{-2em}-\widetilde{\mathbb{O}}^\dagger_3\widetilde{\mathbb{L}}_3=\left[
\begin{array}{c}
    0.4739\\
    0.7689\\
    2.7381
\end{array}
\right],~
-\widetilde{\mathbb{O}}^\dagger_2\widetilde{\mathbb{L}}_2=\left[
\begin{array}{c}
    0.0333\\
    0.0540\\
    0.1923
\end{array}
\right],~\\
&&\hspace{-2em}-\widetilde{\mathbb{O}}^\dagger_1\widetilde{\mathbb{L}}_1=\left[
\begin{array}{c}
    0.0168\\
    0.0273\\
    0.0971
\end{array}
\right],~
-\widetilde{\mathbb{O}}^\dagger_0\widetilde{\mathbb{L}}_0=\left[
\begin{array}{c}
    0.0077\\
    0.0124\\
    0.0443
\end{array}
\right],\\
&&\hspace{-2em}-\widetilde{\mathbb{O}}^\dagger_3\widetilde{\theta}_3=\left[
\begin{array}{c}
    0.4739\\
    0.7689\\
    2.7381
\end{array}
\right],~
-\widetilde{\mathbb{O}}^\dagger_2\widetilde{\theta}_2=\left[
\begin{array}{c}
    0.1221\\
    0.1981\\
    0.7055
\end{array}
\right],~\\
&&\hspace{-2em}-\widetilde{\mathbb{O}}^\dagger_1\widetilde{\theta}_1=\left[
\begin{array}{c}
    0.0922\\
    0.1496\\
    0.5328
\end{array}
\right],~
-\widetilde{\mathbb{O}}^\dagger_0\widetilde{\theta}_0=\left[
\begin{array}{c}
    0.0730\\
    0.1185\\
    0.4220
\end{array}
\right],
\end{eqnarray*}
Then, the unique feedback equilibrium strategy is given by $\{(-\widetilde{\mathbb{O}}^\dagger_k \widetilde{\mathbb{L}}_k, -\widetilde{\mathbb{O}}^\dagger_k \widetilde{\theta}_k), k=0,1,2,3\}$.

\textbf{\underline{Mixed equilibrium solution}}

We use the command ``randn" of MATLAB to randomly generate a $\Phi=\{\Phi_k, k=0,1,2,3\}$. Notie that $\Phi_k\in \mathbb{R}^{3\times 1}$, let $\phi=(\Phi_0,\Phi_1, \Phi_2,\Phi_3)$. We perform the iterations (\ref{S-MV-mixed})-(\ref{pi-MV-mixed}) for 10 times, and get the following are the 10 $\phi$'s and the eigenvalues of the corresponding $\mathcal{O}_k, k=0,1,2,3$,
\begin{eqnarray*}
&&\hspace{-1.5em}\left\{
\begin{array}{l}\phi=\left[
\begin{array}{cccc}
    0.0335&   -1.5771&   -2.0518 &   1.0984\\
   -1.3337&    0.5080&   -0.3538 &  -0.2779\\
    1.1275&    0.2820&   -0.8236 &   0.7015
\end{array}
\right],~~\\
\left\{
\begin{array}{ll}
    0.0099,~   0.0712, ~   0.2092:& \mbox{eigenvalues of }\mathcal{O}_0\\
    0.1309,~    0.0062,~    0.0446:&\mbox{eigenvalues of }\mathcal{O}_1\\
    0.0094,~    0.0759,~    0.2228:&\mbox{eigenvalues of }\mathcal{O}_2\\
    0.0041,~    0.0318,~    0.0930:&\mbox{eigenvalues of }\mathcal{O}_3\\
\end{array}
\right.
\end{array}
\right.\\
&&\hspace{-1.5em}\left\{
\begin{array}{l}\phi=\left[
\begin{array}{cccc}
   -0.5336&   -0.8314&   -0.2620&    0.3502\\
   -2.0026&   -0.9792&   -1.7502&   -0.2991\\
    0.9642&   -1.1564&   -0.2857&    0.0229
\end{array}
\right],~~\\
\left\{
\begin{array}{ll}
    0.0549,~    0.0040, ~  -0.0004:&~~~~\mathcal{O}_0,\\
    0.0054,~   0.0361,  ~  0.1087:&~~~~\mathcal{O}_1,\\
    0.0073,~    0.0574, ~   0.1677:&~~~~\mathcal{O}_2,\\
    0.0041,~    0.0318,~    0.0930:&~~~~\mathcal{O}_3,
\end{array}
\right.
\end{array}
\right.
\\
&&\hspace{-1.5em}\left\{
\begin{array}{l}\phi=\left[
\begin{array}{cccc}
   -0.5890&   -0.7145&   -0.7982&    0.5201\\
   -0.2938&    1.3514&    1.0187&   -0.0200\\
   -0.8479&   -0.2248&   -0.1332&   -0.0348
\end{array}
\right],~~\\
\left\{
\begin{array}{ll}
    0.0199, ~   0.1601,~    0.4696:&~~~~\mathcal{O}_0,\\
    0.0129, ~   0.1032,~    0.3022:&~~~~\mathcal{O}_1,\\
    0.0080, ~   0.0630,~    0.1843:&~~~~\mathcal{O}_2,\\
    0.0041, ~   0.0318,~   0.0930:&~~~~\mathcal{O}_3,
\end{array}
\right.
\end{array}
\right.\\
&&\hspace{-1.5em}\left\{
\begin{array}{l}\phi=\left[
\begin{array}{cccc}
   -2.3299 &  -0.1765&    0.3075&   -1.1201\\
   -1.4491 &   0.7914&   -1.2571&    2.5260\\
    0.3335 &  -1.3320&   -0.8655&    1.6555
\end{array}
\right],~~\\
\left\{
\begin{array}{ll}
    0.0124, ~   0.0986,~    0.2887:&~~~~\mathcal{O}_0,\\
    0.0105, ~   0.0847,~    0.2484:&~~~~\mathcal{O}_1,\\
    0.0139, ~   0.1157, ~   0.3434:&~~~~\mathcal{O}_2,\\
    0.0041, ~   0.0318, ~   0.0930:&~~~~\mathcal{O}_3,
\end{array}
\right.
\end{array}
\right.\\
&&\hspace{-1.5em}\left\{
\begin{array}{l}\phi=\left[
\begin{array}{cccc}
   -0.3349 &  -1.3617&    0.1837 &   0.3914\\
    0.5528 &   0.4550&   -0.4762 &   0.4517\\
    1.0391 &  -0.8487&    0.8620 &  -0.1303
\end{array}
\right],~~\\
\left\{
\begin{array}{ll}
    0.0132,~    0.1042, ~   0.3045:&~~~~\mathcal{O}_0,\\
    0.0135,~    0.1088, ~   0.3192:&~~~~\mathcal{O}_1,\\
    0.0086,~    0.0683, ~   0.2000:&~~~~\mathcal{O}_2,\\
    0.0041,~    0.0318, ~   0.0930:&~~~~\mathcal{O}_3,
\end{array}
\right.
\end{array}
\right.\\
&&\hspace{-1.5em}\left\{
\begin{array}{l}\phi=\left[
\begin{array}{cccc}
    0.8261&   -0.3031&   -0.0679&   -1.1176\\
    1.5270&    0.0230&   -0.1952&    1.2607\\
    0.4669&    0.0513&   -0.2176&    0.6601
\end{array}
\right],~~\\
\left\{
\begin{array}{ll}
    0.0157, ~   0.1258,~    0.3686:&~~~~\mathcal{O}_0,\\
    0.0122, ~   0.0979,~    0.2870:&~~~~\mathcal{O}_1,\\
    0.2324, ~   0.0098,~    0.0791:&~~~~\mathcal{O}_2,\\
    0.0041, ~   0.0318,~    0.0930:&~~~~\mathcal{O}_3,
\end{array}
\right.
\end{array}
\right.\\
&&\hspace{-1.5em}\left\{
\begin{array}{l}\phi=\left[
\begin{array}{cccc}
   -0.9415&    0.1352 &  -1.0298 &  -0.2097\\
   -0.1623&    0.5152 &   0.9492 &   0.6252\\
   -0.1461&    0.2614 &   0.3071 &   0.1832
\end{array}
\right],~~\\
\left\{
\begin{array}{ll}
    0.0218,~    0.1776, ~   0.5230:&~~~~\mathcal{O}_0,\\
    0.0143,~    0.1158, ~   0.3401:&~~~~\mathcal{O}_1,\\
    0.0088,~   0.0701,  ~  0.2053:&~~~~\mathcal{O}_2,\\
    0.0041,~    0.0318, ~   0.0930:&~~~~\mathcal{O}_3,
\end{array}
\right.
\end{array}
\right.\\
&&\hspace{-1.5em}\left\{
\begin{array}{l}\phi=\left[
\begin{array}{cccc}
   -0.2490 &  -0.1922&   -0.4838&   -0.5320\\
   -1.0642 &  -0.2741&   -0.7120&    1.6821\\
    1.6035 &   1.5301&   -1.1742&   -0.8757
\end{array}
\right],~~\\
\left\{
\begin{array}{ll}
    0.0134,~    0.1005, ~   0.2933:&~~~~\mathcal{O}_0,\\
    0.0064,~    0.0463, ~   0.1358:&~~~~\mathcal{O}_1,\\
    0.0086,~    0.0684, ~   0.2001:&~~~~\mathcal{O}_2,\\
    0.0041,~    0.0318, ~   0.0930:&~~~~\mathcal{O}_3,
\end{array}
\right.\end{array}
\right.\\
&&\hspace{-1.5em}\left\{
\begin{array}{l}\phi=\left[
\begin{array}{cccc}
   -1.2507&   -0.2612&   -0.4446&    1.2347\\
   -0.9480&    0.4434&   -0.1559&   -0.2296\\
   -0.7411&    0.3919&    0.2761&   -1.5062
\end{array}
\right],~~\\
\left\{
\begin{array}{ll}
    0.0166, ~   0.1257, ~   0.3667:&~~~~\mathcal{O}_0,\\
    0.0094, ~   0.0701, ~   0.2049:&~~~~\mathcal{O}_1,\\
    0.0058, ~   0.0435, ~   0.1271:&~~~~\mathcal{O}_2,\\
    0.0041, ~   0.0318, ~   0.0930:&~~~~\mathcal{O}_3,
\end{array}
\right.
\end{array}
\right.\\
&&\hspace{-1.5em}\left\{
\begin{array}{l}\phi=\left[
\begin{array}{cccc}
   -0.0290&   -1.0667&   -3.0292&   -0.5078\\
    0.1825&    0.9337&   -0.4570&   -0.3206\\
   -1.5651&    0.3503&    1.2424&    0.0125
\end{array}
\right],~~\\
\left\{
\begin{array}{ll}
    0.0143, ~   0.1050, ~   0.3073:&~~~~\mathcal{O}_0,\\
    0.0075, ~   0.0534, ~   0.1571:&~~~~\mathcal{O}_1,\\
    0.0063, ~   0.0481, ~   0.1404:&~~~~\mathcal{O}_2,\\
    0.0041, ~   0.0318, ~   0.0930:&~~~~\mathcal{O}_3.
\end{array}
\right.
\end{array}
\right.
\end{eqnarray*}
For each $\phi$, all the eigenvalues of $\mathcal{O}_k, k=0,1,2,3$, are nonzero; this means that $\mathcal{O}_k, k=0,1,2,3$, are all invertible. Therefore,
the corresponding (\ref{T-MV-mixed}) (\ref{pi-MV-mixed}) are solvable. In addition, (\ref{S-MV-mixed}) is clearly solvable. From Proposition \ref{Proposition-MV-mixed}, we know that for all the above 10 cases the mixed equilibrium solutions exist.

For example, with the last $\phi$ above, the  mixed equilibrium solution is as follows. Let
\begin{eqnarray*}
&&\Phi_0=\left[
\begin{array}{c}
-0.0290\\    0.1825\\   -1.5651
\end{array}
\right],~~
\Phi_1=\left[
\begin{array}{c}
-1.0667\\0.9337\\0.3503
\end{array}
\right],~~\\
&&\Phi_2=\left[
\begin{array}{c}
-3.0292\\-0.4570\\1.2424
\end{array}
\right],~~
\Phi_3=\left[
\begin{array}{c}
 -0.5078\\-0.3206\\0.0125
\end{array}
\right],
\end{eqnarray*}
and
$$v^{0,x}_k=-({\mathcal{O}}_k^\dagger {\mathcal{L}}_k+\Phi_k)X^{0,x,*}_k-{\mathcal{O}}_k^\dagger {\theta}_k$$
 with
$$
\left\{
\begin{array}{l}
{X}^{0,x,*}_{k+1} =\big{(}s_k-O_k^T{\mathcal{O}}_k^{\dagger}{\mathcal{L}}_k \big{)}{X}^{0,x,*}_{k}-O_k^T{\mathcal{O}}^\dagger_k{\theta}_k,\\[1mm]
{X}^{0,x,*}_{0} = x, k=0, 1, 2, 3,
\end{array}
\right.$$ and
\begin{eqnarray*}
&&\hspace{-1.5em}-\mathcal{O}_0^\dagger\mathcal{L}_0=\left[
\begin{array}{c}
    0.2274\\
    0.3689\\
    1.3137
\end{array}
\right],~~
-\mathcal{O}_1^\dagger\mathcal{L}_1=\left[
\begin{array}{c}
    0.3611\\
    0.5858\\
    2.0862
\end{array}
\right],~~\\
&&\hspace{-1.5em}-\mathcal{O}_2^\dagger\mathcal{L}_2=\left[
\begin{array}{c}
    0.3382\\
    0.5486\\
    1.9537
\end{array}
\right],~~
-\mathcal{O}_3^\dagger\mathcal{L}_3=\left[
\begin{array}{c}
    0.4739\\
    0.7689\\
    2.7381
\end{array}
\right],\\
&&\hspace{-1.5em}-\mathcal{O}_0^\dagger\theta_0=\left[
\begin{array}{c}
    0.2195\\
    0.3561\\
    1.2683
\end{array}
\right],~~
-\mathcal{O}_1^\dagger\theta_1=\left[
\begin{array}{c}
    0.3543\\
    0.5747\\
    2.0468
\end{array}
\right],~~\\
&&\hspace{-1.5em}-\mathcal{O}_2^\dagger\theta_2=\left[
\begin{array}{c}
    0.3365\\
    0.5460\\
    1.9443
\end{array}
\right],~~
-\mathcal{O}_3^\dagger\theta_3=\left[
\begin{array}{c}
    0.4739\\
    0.7689\\
    2.7381
\end{array}
\right].
\end{eqnarray*}
Then, $(\Phi, v^{0,x})$ is a mixed equilibrium portfolio solution of Problem (MV) for $(0,x)$, where $\Phi=\{\Phi_k, k=0,1,2,3\}$.

%
%
%
%
%



%

\end{document}